\documentclass[%
 reprint,
 amsmath,amssymb,
 aps,
 pra,
superscriptaddress]{revtex4-2}

\UseRawInputEncoding

\usepackage[colorlinks=true, allcolors=blue]{hyperref}
\usepackage{amsmath}
\usepackage{braket}
\usepackage{graphicx}
\usepackage{dcolumn}
\usepackage{bm}
\usepackage{listings}
\usepackage{xcolor}

\begin{document}
\preprint{APS/123-QED}
\title{Numerical Exploration of the Pythagorean Theorem Using HOBO Algorithm}

\author{Shoya Yasuda}
\email{yasuda@vigne-cla.com}
\affiliation{Vignette \& Clarity, Inc., 6-24-2 Honkomagome, Bunkyo-ku, Tokyo, 113-0021, Japan}
\affiliation{Tokyo Institute of Technology, School of Computing, J2 bldg. room1710, 4259 Nagatsuta-cho, Midori-ku, Yokohama, 226-8503, Japan}
\author{Naoaki Mochida}
\affiliation{St. Mary's International School, 1-6-19 Seta Setagaya-ku, Tokyo, 158-8668, Japan}
\author{Shunsuke Sotobayashi}
\email{derwind0707@gmail.com}
\affiliation{Independent Researcher}
\author{Devanshu Garg}
\email{devanshu@blueqat.com}
\affiliation{blueqat Inc., 2-24-12 Shibuya, Shibuya-ku, Tokyo, 150-6139, Japan}
\author{Yuichiro Minato}
\email{minato@blueqat.com}
\affiliation{blueqat Inc., 2-24-12 Shibuya, Shibuya-ku, Tokyo, 150-6139, Japan}
\begin{abstract}
This paper introduces a novel method for finding integer sets that satisfy the Pythagorean theorem by leveraging the Higher-Order Binary Optimization (HOBO) formulation. Unlike the Quadratic Unconstrained Binary Optimization (QUBO) formulation, which struggles to express complex mathematical equations, HOBO's ability to model higher-order interactions between binary variables makes it well-suited for addressing more complex and expressive problem settings. \end{abstract}
\maketitle
\section{Introduction}
In recent years, the application of quantum computing and classical solvers has advanced, leading to the exploration of new approaches for efficiently solving optimization problems and mathematical challenges. Among these, the QUBO (Quadratic Unconstrained Binary Optimization) formulation, which optimizes quadratic functions of binary variables~\cite{Andrew2014}, has been widely adopted in solving a variety of mathematical problems.

However, with the ongoing advancements in quantum computing technology, the HOBO (Higher-Order Binary Optimization) formulation is gaining recognition for its enhanced problem-solving capabilities. HOBO allows for the optimization of higher-order binary functions (third-order and beyond), making it well-suited for tackling more complex and expressive problem settings. Particularly, HOBO holds promise for new developments, including integer encoding.

In this study, rather than focusing on societal issues~\cite{Domino_2022}~\cite{Glos2022}~\cite{Tabi_2020}, we investigate the utility of the HOBO formulation for solving general mathematical problems. Specifically, we explored the feasibility and effectiveness of implementing HOBO using classical HOBO solvers~\cite{minato2024tensornetworkbasedhobo}~\cite{yasuda2024hobotanefficienthigherorder}. Through this investigation, we demonstrated that HOBO is not merely a theoretical model but also an effective tool for addressing real-world mathematical problem-solving. The HOBO formulation, with its expressive power and flexibility, offers a new approach to mathematical problem-solving and holds the potential to complement and enhance traditional QUBO-based methods.

\section{HOBO}
Higher-Order Binary Optimization (HOBO) is a mathematical framework designed to solve higher-order binary optimization problems. Unlike QUBO (Quadratic Unconstrained Binary Optimization), which is limited to quadratic terms involving binary variables, HOBO extends this to include third-order and higher interactions, enabling the modeling of more complex problem settings.

\subsection{HOBO Formulation}
The general form of a HOBO problem can be expressed as the problem of finding \( x \) that minimizes \( f(x) \):

\begin{equation}
f(\mathbf{x}) = \sum_i c_i x_i + \sum_{ij} c_{ij} x_i x_j + \sum_{ijk} c_{ijk} x_i x_j x_k + \dots 
\end{equation}

Here, 
\begin{itemize}
    \item $\mathbf{x} = \{x_1, x_2, \dots, x_n\}$ is a vector of binary variables where each $x_i \in \{0,1\}$.
    \item $c_i$, $c_{ij}$, $c_{ijk}$, etc., are coefficients representing the weights of the terms.
    \item The expression includes linear, quadratic, cubic, and potentially higher-order interaction terms.
\end{itemize}

\subsection{Example: Third-Order HOBO}
As an example, consider a third-order HOBO problem:

\begin{equation}
\quad f(\mathbf{x}) = c_1 x_1 + c_{12} x_1 x_2 + c_{123} x_1 x_2 x_3
\end{equation}

In this case, $f(\mathbf{x})$ includes the cubic term $c_{123} x_1 x_2 x_3$, which qualifies it as a HOBO problem. This higher-order term is something that cannot be handled by traditional QUBO formulations without decomposition.

\section{Methods for Solving HOBO Problems}
There are some approaches to solving HOBO (Higher-Order Binary Optimization) problems: utilizing the Quantum Approximate Optimization Algorithm (QAOA) by transforming the cost Hamiltonian from HOBO, and employing classical solvers designed to handle HOBO directly.

\subsection{QAOA Approach}
The Quantum Approximate Optimization Algorithm (QAOA) is a quantum algorithm that is well-suited for solving combinatorial optimization problems. To apply QAOA to a HOBO problem, the cost Hamiltonian, which represents the objective function, must first be transformed from the HOBO formulation into a form compatible with QAOA.

\subsection{Classical Solvers for HOBO}
Alternatively, HOBO problems can be addressed using classical solvers that are specifically designed or adapted to handle higher-order binary optimization.
There are classical algorithms that can directly solve HOBO problems by efficiently navigating the search space of higher-order interactions. 
Classical solvers offer a practical alternative to quantum algorithms, for problems where quantum resources are limited or where the problem size is within the reach of classical computational power.

\section{TYTAN (HOBOTAN) Classical Solver}
In this work, we introduce TYTAN (HOBOTAN) [ https://pypi.org/project/tytan/ ], a classical solver specifically developed to address higher-order binary optimization (HOBO) problems. By utilizing tensor networks~\cite{Vidal_2003}, TYTAN is capable of efficiently solving complex higher-order problems.

\subsection{Compilation of HOBO into Tensor Networks}
Similar to the way QUBO problems are handled, HOBO problems are formulated and then compiled into a tensor representation within the TYTAN solver. The key process involves transforming the HOBO problem into a tensor format, where the higher-order interactions among binary variables are encoded as tensors.

\begin{itemize}
    \item \textbf{HOBO Formulation:} The HOBO problem is first formulated with the necessary higher-order terms, which are then mapped onto a corresponding tensor network. This transformation preserves the problem's structure while making it suitable for tensor-based computation.
    \item \textbf{Tensor Compilation:} Once the HOBO problem is formulated, the TYTAN solver compiles it into a tensor network~\cite{wang2023tensornetworksmeetneural}. This step is critical as it allows the solver to efficiently represent and manipulate the higher-order interactions present in the problem.
\end{itemize}

\subsection{Contraction Calculations}
Within the TYTAN solver, the actual computation is performed through contraction calculations between the problem tensor and the solution vector. This approach allows TYTAN to handle high-dimensional problems effectively.

\begin{itemize}
    \item \textbf{Tensor Contraction:} The solver performs contraction operations between the tensors representing the HOBO problem and the vectors that represent possible solutions. By systematically contracting the tensors, TYTAN can explore the solution space efficiently and identify optimal solutions.
    \item \textbf{Handling Higher-Order Problems:} This contraction process is particularly well-suited for higher-order problems, as it allows the solver to manage and compute the complex interactions encoded in the tensor network, providing a robust solution to HOBO challenges.
\end{itemize}

\section{Advantages of TYTAN (HOBOTAN)}
TYTAN leverages the power of tensor networks to efficiently solve higher-order optimization problems that are beyond the capabilities of traditional methods. By compiling HOBO problems into tensor networks and performing contraction calculations, TYTAN provides a powerful tool for tackling complex optimization challenges in a classical computing environment.

\section{Pythagorean Theorem}
The Pythagorean Theorem is a fundamental principle in Euclidean geometry, which establishes a relationship between the lengths of the sides of a right-angled triangle. The theorem is named after the ancient Greek mathematician Pythagoras and is one of the most well-known and widely used mathematical theorems.

\subsection{Statement of the Theorem}
The Pythagorean Theorem states that in a right-angled triangle, the square of the length of the hypotenuse (the side opposite the right angle) is equal to the sum of the squares of the lengths of the other two sides. Mathematically, the theorem is expressed as:

\begin{equation}
z^2 = x^2 + y^2
\end{equation}

where:
\begin{itemize}
    \item $z$ is the length of the hypotenuse.
    \item $x$ and $y$ are the lengths of the other two sides.
\end{itemize}

\subsection{Example}
Consider a right-angled triangle where the lengths of the two shorter sides are 3 and 4. According to the Pythagorean Theorem, the length of the hypotenuse can be calculated as follows:

\begin{equation}
z^2 = 3^2 + 4^2 = 9 + 16 = 25
\end{equation}

Taking the square root of both sides:

\begin{equation}
z = \sqrt{25} = 5
\end{equation}

Thus, the hypotenuse of the triangle is 5.

\section{Steps to Find Numbers Satisfying the Pythagorean Theorem Using HOBO}
This section outlines the steps to find numbers that satisfy the Pythagorean theorem using the Higher-Order Binary Optimization (HOBO) framework.

\subsection{Step 1: Convert Binary HOBO to Integers Using Integer Encoding}
First, to utilize HOBO, we need to encode the integer values into binary form. This encoding process involves representing the integer variables \(x\), \(y\), and \(z\) as binary variables.

For example, if we represent each variable using 4 bits, the integers \(x\), \(y\), and \(z\) can be expressed as:

\begin{equation}
\begin{aligned}
x &= 2^0x_0 + 2^1x_1 + 2^2x_2 + 2^3x_3 \\
y &= 2^0y_0 + 2^1y_1 + 2^2y_2 + 2^3y_3 \\
z &= 2^0z_0 + 2^1z_1 + 2^2z_2 + 2^3z_3
\end{aligned}
\end{equation}

Here, \(x_i\), \(y_i\), and \(z_i\) are binary variables, and the integers are derived through integer encoding.

\subsection{Step 2: Transform the Pythagorean Theorem into a Minimization Problem}
Next, to find integer values that satisfy the Pythagorean theorem, we construct the following equation:

\begin{equation}
H = (x^2 + y^2 - z^2)^2
\end{equation}

This equation is structured to yield a value of zero when the Pythagorean theorem \(x^2 + y^2 = z^2\) is satisfied. Minimizing this equation corresponds to finding the integer values \(x\), \(y\), and \(z\) that satisfy the Pythagorean theorem.

\subsection{Step 3: Expand the Equation and Convert to HOBO Form}
The next step involves expanding the equation \(H\) and converting it into HOBO form. The HOBO form is a representation that includes higher-order interaction terms among binary variables.

The expanded equation takes the following form:

\begin{equation}
H = (x^2 + y^2 - z^2)^2 = x^4 + 2x^2y^2 + y^4 - 2x^2z^2 - 2y^2z^2 + z^4
\end{equation}

This equation includes higher-order interaction terms (e.g., \(x^4\) and \(x^2y^2\)). Converting this expanded equation into HOBO form prepares it for processing by a HOBO solver.

\subsection{Step 4: Input the Formulation into a HOBO Solver to Find the Solution}
Finally, the equation in HOBO form is input into a HOBO solver. The solver will minimize the equation, thereby finding the integer values \(x\), \(y\), and \(z\) that satisfy the Pythagorean theorem.

For example, using a solver like TYTAN (HOBOTAN), the solution can be found by sampling for the lowest energy state, which corresponds to the set of integers that satisfy the Pythagorean theorem.

By following these steps, the HOBO framework can effectively be used to search for integers that satisfy the Pythagorean theorem.

\section{Numerical Solution and Results}
In this section, we detail the results obtained from solving the Pythagorean theorem problem using the HOBO framework.

\subsection{Example: Reformulation of Variables}

In the original formulation, the variables \(x\), \(y\), and \(z\) were represented as shown in equation (6). However, to ensure that \(x\), \(y\), and \(z\) are all greater than or equal to 1, we reformulated the variables by subtracting 1 from each expression. This allowed us to reduce the number of quantum bits required for the encoding by 1, resulting in the following expressions:

\begin{equation}
\begin{aligned}
x &= 1 + (2^0x_0 + 2^1x_1 + 2^2x_2) \\
y &= 1 + (2^0y_0 + 2^1y_1 + 2^2y_2) \\
z &= 1 + (2^0z_0 + 2^1z_1 + 2^2z_2)
\end{aligned}
\end{equation}

In this revised formulation, \(x\), \(y\), and \(z\) are guaranteed to be at least 1. It eliminates the need for additional constraints to ensure that the values are greater than or equal to 1, providing a significant advantage.

\subsection{Precision Issues in PyTorch Compared to Numpy}

Typically, Numpy uses ``float64'' as the default floating-point type, while PyTorch uses ``float32'' by default. This difference led to a situation where calculations that returned correct values in Numpy produced incorrect results when run in PyTorch's GPU mode.

In our experiment to find Pythagorean numbers, where \( x \), \( y \), and \( z \) are each defined using 7 bits or more, we found that some solutions were missing. This issue seems to arise during energy calculations, though it is likely that the HOBO tensors themselves are functioning correctly. The unexpected behavior is believed to be due to the einsum operations.

It is generally known that the float type can accurately represent integers up to 7 digits, while the double type can represent integers up to 15 digits accurately. In our case, calculations involving values around \( 2^7 \) began to yield inaccurate results, highlighting the importance of using ``float64'' for accurate computation.

\subsection{Problem Setting and Solution:1}

The goal is to find all Pythagorean triples \((x, y, z)\) where each of \(x\), \(y\), and \(z\) are less than 16. The expected solutions are the four Pythagorean triples: \((3, 4, 5)\), \((6, 8, 10)\), \((5, 12, 13)\), and \((9, 12, 15)\). Considering that \(x\) and \(y\) can be swapped, there are a total of 8 possible solutions.

To solve this problem, we encoded \(x\), \(y\), and \(z\) and applied the HOBO framework to efficiently identify all valid Pythagorean triples within the specified range. The solution process confirmed the expected 8 solutions, covering all permutations of the 4 unique triples.

The Python code that was actually executed is as follows:

\begin{lstlisting}[language=Python]
from tytan import *

# Define quantum bits
qx = symbols_list(4, 'qx{}')
qy = symbols_list(4, 'qy{}')
qz = symbols_list(4, 'qz{}')

# Represent x, y, z in binary form
x = 1 + 1*qx[0] + 2*qx[1] + 4*qx[2] + 8*qx[3]
y = 1 + 1*qy[0] + 2*qy[1] + 4*qy[2] + 8*qy[3]
z = 1 + 1*qz[0] + 2*qz[1] + 4*qz[2] + 8*qz[3]
print(x)

# Pythagorean equation constraint
H = (x**2 + y**2 - z**2)**2

# Compile the HOBO tensor
hobo, offset = Compile(H).get_hobo()
print(f'offset\n{offset}')

# Select the sampler
solver = sampler.MIKASAmpler()

# Sampling
result = solver.run(hobo, shots=10000)

# Display top 10 results
for r in result[:10]:
    print(f'Energy {r[1]}, Occurrence {r[2]}')

    print('x =', Auto_array(r[0]).get_nbit_value(x))
    print('y =', Auto_array(r[0]).get_nbit_value(y))
    print('z =', Auto_array(r[0]).get_nbit_value(z))

\end{lstlisting}

The calculation was performed 10,000 times, and the solutions obtained, which exhibited some variability as part of the sampling process, are shown below.

\begin{verbatim}
1 + qx0 + 2*qx1 + 4*qx2 + 8*qx3
offset
1.0
MODE: GPU
DEVICE: cuda:0
Energy -1.0, Occurrence 1234
x = 8.0
y = 6.0
z = 10.0
Energy -1.0, Occurrence 929
x = 12.0
y = 5.0
z = 13.0
Energy -1.0, Occurrence 775
x = 12.0
y = 9.0
z = 15.0
Energy -1.0, Occurrence 983
x = 4.0
y = 3.0
z = 5.0
Energy -1.0, Occurrence 1021
x = 9.0
y = 12.0
z = 15.0
Energy -1.0, Occurrence 1316
x = 5.0
y = 12.0
z = 13.0
Energy -1.0, Occurrence 1332
x = 3.0
y = 4.0
z = 5.0
Energy -1.0, Occurrence 849
x = 6.0
y = 8.0
z = 10.0
Energy 0.0, Occurrence 223
x = 4.0
y = 7.0
z = 8.0
Energy 0.0, Occurrence 39
x = 14.0
y = 1.0
z = 14.0
\end{verbatim}

Since the constant is 1, the correct solution corresponds to a calculated cost function of \(-1\) after subtracting the constant. In this case, by excluding the solutions with a cost of \(0\), which appeared less frequently, the 8 correct solutions with a cost of \(-1\) were obtained with a high probability as shown in Table\ref{tab:quantum_optimization_results}.

\begin{table}[h]
\centering
\begin{tabular}{|c|c|c|c|c|c|}
\hline
\textbf{\#} & \textbf{x} & \textbf{y} & \textbf{z} & \textbf{Energy} & \textbf{Occurrence} \\ \hline
1 & 8 & 6 & 10 & -1.0 & 1234 \\ \hline
2 & 12 & 5 & 13 & -1.0 & 929 \\ \hline
3 & 12 & 9 & 15 & -1.0 & 775 \\ \hline
4 & 4 & 3 & 5 & -1.0 & 983 \\ \hline
5 & 9 & 12 & 15 & -1.0 & 1021 \\ \hline
6 & 5 & 12 & 13 & -1.0 & 1316 \\ \hline
7 & 3 & 4 & 5 & -1.0 & 1332 \\ \hline
8 & 6 & 8 & 10 & -1.0 & 849 \\ \hline
9 & 4 & 7 & 8 & 0.0 & 223 \\ \hline
10 & 14 & 1 & 14 & 0.0 & 39 \\ \hline
\end{tabular}
\caption{Results of the Optimization Using 4 Qubits}
\label{tab:quantum_optimization_results}
\end{table}

\subsection{Initial Validation with HOBO Model and Further Comparison with QUBO}

We first validated the Pythagorean condition for \((x, y, z)\) sets using the HOBO model. Once the HOBO model successfully confirmed the expected results, we proceeded to compare the performance of the QUBO and HOBO models with larger numbers.

Encouraged by these results, we extended our investigation to larger values and performed a detailed comparison between the QUBO and HOBO models. In the QUBO model, the squares of \(x\), \(y\), and \(z\) were prepared, and quantum bits corresponding to each of these numbers were defined. A strong one-hot constraint was applied, along with a weak equation constraint of \(x^2 + y^2 = z^2\). However, it is important to note that in the QUBO model, \(x\), \(y\), and \(z\) cannot be directly defined with quantum bits due to the fourth-order terms produced by the equation constraint \((7)\).

The Python code that was actually executed for the QUBO model is as follows:

\begin{lstlisting}[language=Python]
import numpy as np
from tytan import *

Power = 4
x = np.arange(1, 1 + 2**Power)
y = np.arange(1, 1 + 2**Power)
z = np.arange(1, 1 + 2**Power)

x2 = x**2
y2 = y**2
z2 = z**2

# Define quantum bits
qx = symbols_list(len(x2), 'qx{}')
qy = symbols_list(len(y2), 'qy{}')
qz = symbols_list(len(z2), 'qz{}')

H = (sum(qx) - 1)**2
H += (sum(qy) - 1)**2
H += (sum(qz) - 1)**2

H += 0.01 * (sum(x2*qx) + sum(y2*qy) - sum(z2*qz))**2

# Compile the QUBO matrix
qubo, offset = Compile(H).get_qubo()
print(f'offset\n{offset}')

# Select the sampler
solver = sampler.ArminSampler()

# Sampling
result = solver.run(qubo, shots=10000)

# Display top 10 results
ans = []
for r in result[:10]:
    print(f'Energy {r[1]}, Occurrence {r[2]}')
    
    #check one-hot constraint
    xs, xsubs = Auto_array(r[0]).get_ndarray('qx{}')
    ys, ysubs = Auto_array(r[0]).get_ndarray('qy{}')
    zs, zsubs = Auto_array(r[0]).get_ndarray('qz{}')
    
    if sum(xs) != 1:
        print('ERR: x is not One-hot.')
        continue
    if sum(ys) != 1:
        print('ERR: y is not One-hot.')
        continue
    if sum(zs) != 1:
        print('ERR: z is not One-hot.')
        continue
    
    print('x =', x[np.argmax(xs)])
    print('y =', y[np.argmax(ys)])
    print('z =', z[np.argmax(zs)])
\end{lstlisting}

As part of this extended comparison, we systematically increased the bit length of \(x\), \(y\), and \(z\), and analyzed the number of quantum bits required for each model. The Power of two in our analysis refers to the maximum value that \(x\), \(y\), and \(z\) can represent as powers of two. For example, Power = 4 indicates that the models explore all integers within the range \(1 \leq x \leq 16\), \(1 \leq y \leq 16\), and \(1 \leq z \leq 16\). The number of quantum bits required by the QUBO and HOBO models to represent integers as powers of two is shown in Table \ref{tab:qubit_comparison}.

\begin{table}[ht]
\centering
\begin{tabular}{|c|c|c|}
\hline
\textbf{Power of two} & \textbf{QUBO} & \textbf{HOBO} \\ \hline
3                     & 24            & 9             \\ \hline
4                     & 48            & 12            \\ \hline
5                     & 96            & 15            \\ \hline
6                     & 192           & 18            \\ \hline
7                     & 384           & 21            \\ \hline
8                     & 768           & 24            \\ \hline
9                     & 1536          & 27            \\ \hline
10                    & 3072          & 30            \\ \hline
11                    & 6144          & 33            \\ \hline
12                    & 12288         & 36            \\ \hline
\end{tabular}
\caption{Quantum Bits Required}
\label{tab:qubit_comparison}
\end{table}

Through this comparative study, we aimed to evaluate the scalability and accuracy of both the QUBO and HOBO models when applied to larger and more complex problem instances.

The results of exploring primitive Pythagorean triples within the range of \(1 \leq z \leq 16\) (Power = 4) are shown in Table \ref{tab:pythagorean_triples}. The number of shots was set uniformly to 100,000. The Table lists all theoretically possible primitive Pythagorean triples, and the Occ column aggregates the number of times each solution was reached out of 100,000 samplings. It is important to note that for (3, 4, 5), the occurrences of (4, 3, 5) are also included in the count. Additionally, non-primitive Pythagorean triples, such as (6, 8, 10), were discarded and not included in the aggregation.

As observed, the QUBO model only reached the solution a few dozen times out of 100,000 shots even at the Power = 4 stage. In contrast, the HOBO model shows much more robust exploration.

\begin{table}[ht]
\centering
\begin{tabular}{|c|c|c|c|c|}
\hline
\textbf{x} & \textbf{y} & \textbf{z} & \textbf{QUBO-Occ} & \textbf{HOBO-Occ} \\ \hline
3          & 4          & 5          & 32                & 27643             \\ \hline
5          & 12         & 13         & 76                & 19986             \\ \hline
\end{tabular}
\caption{Primitive Pythagorean Triples Discovered in HOBO Model with Power = 4}
\label{tab:pythagorean_triples}
\end{table}

Similarly, the exploration results for Power = 6 are shown in Table \ref{tab:pythagorean_triples_power6}. It can be observed that the HOBO model still retains the ability to discover all primitive Pythagorean triples within the range \(1 \leq z \leq 64\).

\begin{table}[ht]
\centering
\begin{tabular}{|c|c|c|c|c|}
\hline
\textbf{x} & \textbf{y} & \textbf{z} & \textbf{QUBO-Occ} & \textbf{HOBO-Occ} \\ \hline
3          & 4          & 5          & 0                 & 1537             \\ \hline
5          & 12         & 13         & 0                 & 1654             \\ \hline
8          & 15         & 17         & 0                 & 912              \\ \hline
7          & 24         & 25         & 0                 & 703              \\ \hline
20         & 21         & 29         & 0                 & 1474             \\ \hline
12         & 35         & 37         & 0                 & 707              \\ \hline
9          & 40         & 41         & 0                 & 1131             \\ \hline
28         & 45         & 53         & 0                 & 609              \\ \hline
11         & 60         & 61         & 0                 & 552              \\ \hline
\end{tabular}
\caption{Primitive Pythagorean Triples Discovered in HOBO Model with Power = 6}
\label{tab:pythagorean_triples_power6}
\end{table}

\begin{figure}[ht]
\centering
\includegraphics[width=\linewidth]{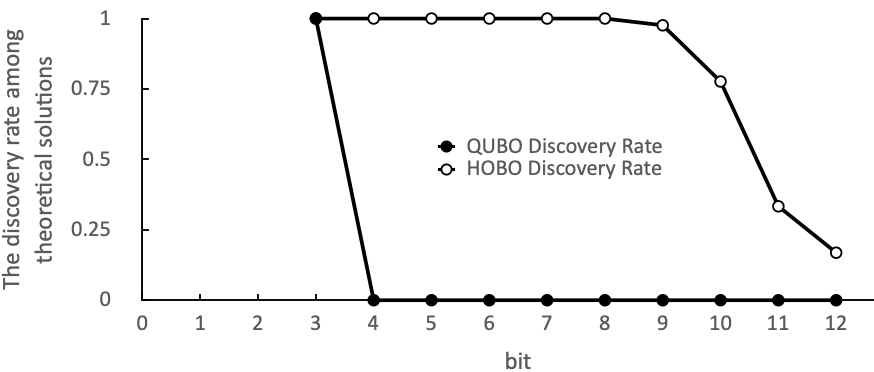} 
\caption{The discovery rate of primitive Pythagorean triples as the Power increases.}
\label{fig:discovery_rate_power}
\end{figure}

The figure \ref{fig:discovery_rate_power} and Table \ref{tab:discovery_rate} show the discovery rate of primitive Pythagorean triples as the Power increases. For example, at Power = 5, where the range is \(1 \leq z \leq 32\), there are theoretically 5 types of primitive Pythagorean triples. The QUBO model failed to find any, whereas the HOBO model successfully discovered all of them.

In the HOBO model, missed discoveries begin to appear from Power = 9, and at Power = 12 (where \(1 \leq z \leq 4096\)), only about 17\% of the primitive Pythagorean triples were discovered. Nonetheless, the results obtained with the HOBO model were significantly better compared to the QUBO model. This improvement is likely due to the simplification of the search space, which was achieved primarily by reducing the number of required quantum bits.

\begin{table}[ht]
\centering
\caption{Discovery Rate of Primitive Pythagorean Triples with Increasing Power}
\begin{tabular}{|c|p{2cm}|p{2cm}|p{2cm}|}
\hline
\textbf{Power of two} & \textbf{Number of Theoretical Solutions} & \textbf{QUBO Discovery Rate} & \textbf{HOBO Discovery Rate} \\ \hline
3                     & 1                                      & 1                            & 1                            \\ \hline
4                     & 2                                      & 0                            & 1                            \\ \hline
5                     & 5                                      & 0                            & 1                            \\ \hline
6                     & 9                                      & 0                            & 1                            \\ \hline
7                     & 20                                     & 0                            & 1                            \\ \hline
8                     & 39                                     & 0                            & 1                            \\ \hline
9                     & 83                                     & 0                            & 0.9759             \\ \hline
10                    & 161                                    & 0                            & 0.7764              \\ \hline
11                    & 327                                    & 0                            & 0.3333            \\ \hline
12                    & 652                                    & 0                            & 0.1687            \\ \hline
\end{tabular}
\label{tab:discovery_rate}
\end{table}

\section{Discussion}
The results obtained demonstrate that the HOBO solver successfully solved the Pythagorean theorem problem within the scope of the HOBO formulation, identifying the correct integer solutions.

The QUBO approach shows some improvement in discovery rate by reducing the number of qubits through leveraging the properties of Pythagorean numbers. However, this improvement is not sufficient to counteract the exponential increase in the number of qubits. Moreover, similar optimizations can be applied to the HOBO model as well. It is important to emphasize that this comparison was conducted using simple implementations.

Additionally, the results may also depend on the specific implementation of the sampling process, so this should be taken into account.

Within the scope of this experiment, no overflow was encountered in the program when using float64.

In this study, we gradually increased the number of qubits used in the encoding while exploring Pythagorean triples with a classical solver. By incrementally increasing the qubit count, we were able to explore a broader range of integers. 

Looking forward, the use of integer encoding within the HOBO framework opens up promising avenues for solving problems that have been challenging to address in the context of quantum computing. By leveraging this method and systematically adjusting the number of qubits, we may overcome some of the limitations encountered in traditional quantum approaches.

\bibliography{apssamp}

\end{document}